\documentclass{article}

\usepackage{delarray,verbatim,enumerate,a4wide}
\usepackage{amsmath,amsthm,amstext,amsbsy,amssymb,amsfonts,amscd}
\usepackage[all]{xy}
\usepackage[frenchb]{babel}
\usepackage[T1]{fontenc}
\usepackage[utf8]{inputenc}
\usepackage{color}
\definecolor{vert}{rgb}{0.1,0.4,0.2}
\usepackage[colorlinks=true,linkcolor=blue,citecolor=vert]{hyperref}

\usepackage{calligra}
\DeclareFontShape{T1}{calligra}{m}{n}{<->s*[0.95]callig15}{}
\DeclareMathAlphabet{\mathscr}{T1}{calligra}{m}{n}

\newtheorem{Th}{Théorème}[]
\newtheorem{Lem}[Th]{Lemme}
\newtheorem{Prop}[Th]{Proposition}
\newtheorem{Cor}[Th]{Corollaire}

\newtheorem{Sco}[Th]{Scolie}

\newtheorem{DTh}[Th]{Définition \& Théorème}

\newtheorem*{Th*}{Théorème} 
\newtheorem*{Sco*}{Scolie}
\newtheorem*{ThC*}{Théorème \& Conjecture}

\def\Preuve{\noindent {\it Preuve.~}}

\def\PreuveCor{\noindent{\it Preuve du Corollaire.~}}

\font\teneufm=eufm10
\font\seveneufm=eufm7
\font\fiveeufm=eufm5
\newfam\gothfam
\textfont\gothfam=\teneufm \scriptfont\gothfam=\seveneufm
\scriptscriptfont\gothfam=\fiveeufm

		\def\QQ{\mathbb Q}	
\def\NN{\mathbb N}	\def\ZZ{\mathbb Z}
\def\F2{\mathbb{F}_2}	\def\Z2{\mathbb{Z}_2}		
\def\Zl{\mathbb{Z}_\ell} 	\def\Ql{\mathbb{Q}_\ell}	

 				\def\U{\mathcal  U}	\def\F{\mathcal  F}
\def\J{\mathcal  J}  	\def\C{\mathcal  C}	\def\R{\mathcal  R}	
 	\def\Pl{\mathcal  P\ell}  	
\def\E{\mathcal  E}				

	\def\p{{\mathfrak p}}	\def\x{{\mathfrak x}}		
		\def\l{{\mathfrak l}}

\def\wi{\widetilde}		
	\def\div{\operatorname{div}}	
	\def\deg{\operatorname{deg}}
\def\Gal{\operatorname{Gal}}	\def\Log{\operatorname{Log}}	
		\def\Hom{\operatorname{Hom}}
\def\Ind{\operatorname{Ind}}

\makeatletter
\newcommand*\wt[2][0.2ex]{%
        \begingroup
        \mathchoice{\wt@helper{#1}{#2}{\displaystyle}{\textfont}}
                   {\wt@helper{#1}{#2}{\textstyle}{\textfont}}
                   {\wt@helper{#1}{#2}{\scriptstyle}{\scriptfont}}
                   {\wt@helper{#1}{#2}{\scriptscriptstyle}{\scriptscriptfont}}%
        \endgroup
        #2%
}
\newcommand*\wt@helper[4]{%
        \def\currentfont{\the#41}%
        \def\currentskewchar{\char\the\skewchar\currentfont}%
        \setbox\tw@\hbox{\currentfont$#2$\currentskewchar}%
        \dimen@ii\wd\tw@
        \setbox\tw@\hbox{\currentfont$#2${}\currentskewchar}%
        \advance\dimen@ii-\wd\tw@
        \rlap{\raisebox{-#1}{$\m@th#3\kern\dimen@ii\widetilde{\phantom{#2}}$}}%
}
\makeatother

\def\wE{\,\wt[0.1ex]{\!\mathcal E}}	\def\we{\wt[0.3ex]{\mathfrak E}}	\def\wU{\wt[0.2ex]{\mathcal U}}
\def\wJ{\,\wt[0.2ex]{\!\mathcal J}}	\def\wCl{\wt[0.1ex]{\mathcal C\!\ell}} \def\wDl{\wt[0.2ex]{\mathcal D\ell}}

\begin{document}

\title{\LARGE\bf Normes cyclotomiques naïves et unités logarithmiques}

\author{ Jean-François {\sc Jaulent} }
\date{}
\maketitle
\bigskip

{\small
\noindent{\bf Résumé.} Nous déterminons le rang du sous-groupe $\wi E_K$ des éléments du groupe multiplicatif d'un corps de nombres $K$ qui sont normes à chaque étage fini de sa $\Zl$-extension cyclotomique $K^c$; et nous comparons son $\ell$-adifié $\Zl\otimes_\ZZ \wi E_K$ avec le $\ell$-groupe des unités logarithmiques $\,\wE_K$. Nous donnons à cette occasion une preuve très facile de la conjecture de Gross-Kuz'min en $\ell$ pour les extensions $K/k$ d'un corps abélien dans lesquelles les places au-dessus de $\ell$ ne se décomposent pas.}\smallskip

{\small
\noindent{\bf Abstract.} We compute the $\ZZ$-rank of the subgroup $\wi E_K =\bigcap_{n\in\NN} N_{K_n/K}(K_n^\times)$ of elements of the multiplicative group of a number field $K$ that are norms from every finite level of the cyclotomic $\Zl$-extension $K^c$ of $K$. Thus we compare its $\ell$-adification  $\Zl\otimes_\ZZ \wi E_K$ with the group of logarithmic units $\,\wE_K$. By the way we point out an easy proof of the Gross-Kuz'min conjecture for $\ell$-undecomposed extensions of abelian fields.}
\bigskip

\noindent{\large \bf Introduction}
\medskip

Le $\ell$-groupe des normes cyclotomiques $\:\wE_K$ attaché à un corps de nombres $K$ se présente comme le sous-module du $\ell$-adifié $\R_K=\Zl\otimes_\ZZ K^\times$ du groupe multiplicatif  $K^\times$ formé des éléments de $\R_K$ qui sont normes dans tous les étages finis $K_n/K$ de la $\Zl$-extension cyclotomique $K^c$ de $K$.\smallskip

Or, comme observé dans \cite{J55}, $\:\wE_K$  ne provient pas en général (par $\ell$-adification) d'un sous-groupe naturel (ou naïf) du groupe $K^\times$\!.
La Théorie $\ell$-adique du corps de classes (cf. \cite{J31} ou \cite{Gra}) permet cependant d'interpréter ce groupe des normes cyclotomiques comme un groupe d'unités, analogue au $\ell$-adifié $\E_K=\Zl\otimes_\ZZ E_K$ du groupe des unités au sens habituel, en lien avec un groupe de classes $\,\wCl_K$, qui est le pendant logarithmique du $\ell$-groupe des classes d'idéaux.\smallskip

Le rang sur $\Zl$ du groupe des unités logarithmiques  $\:\wE_K$ (i.e. la dimension de son quotient par le sous-groupe $\mu_K^{_{(\ell)}}$ des racines de l'unité d'ordre $\ell$-primaire contenues dans $K$) est donné par la formule (cf. \cite{J55}, Sco. 6):

\centerline{$\dim_{\Zl}\wE_K=r_K+c_K+\delta_K^{\mathscr G}$,}\smallskip

\noindent  où  $r_K$ et $c_K$ sont respectivement les nombres de places réelles et complexes de $K$; et $\delta_K^{\mathscr G}$  mesure le défaut dans $K$ de la conjecture de Gross-Kuz'min  pour  le premier $\ell$. En particulier, lorsque le corps $K$ vérifie cette conjecture pour le premier $\ell$  (par exemple pour $K$ abélien sur $\QQ$), il suit:\smallskip

\centerline{ $\wE_K\simeq\mu_K^{_{(\ell)}}  \Zl^{r_K+c_K}$.}\smallskip

Bien entendu, il reste possible, indépendamment de toute conjecture, de définir à l'instar de Bertrandias et Payan (\cite{BP}, \S1.5) le groupe des $\ell$-normes cyclotomiques {\em naïves} comme l'intersection\smallskip

\centerline{$\widetilde E_K =\bigcap_{n\in\NN} N_{K_n/K}(K_n^\times)$}\smallskip

\noindent des sous-goupes normiques de $K^\times$ attachés aux étages finis de la $\Zl$-extension cyclotomique $K^c/K$.\smallskip

Le groupe obtenu est alors un sous-groupe du groupe $E'_K$ des $\ell$-unités du corps $K$, qui contient le $\ZZ$-module multiplicatif engendré par $\ell$, donc un $\ZZ$-module de type fini, produit direct\smallskip

\centerline{$\wi E_K \simeq  \mu_K^{\phantom{l)}}\ZZ^{\tilde e_K}$}\smallskip

\noindent du groupe des racines de l'unité contenues dans $K$ et d'un $\ZZ$-module libre de dimension $\tilde e_K$ avec\smallskip

\centerline{$1 \le \wi e_K =\dim_{\ZZ}\wi E_K \le \dim_{\ZZ}E'_K=r_K+c_K+l_K-1$,}\smallskip

\noindent où $ l_K$ désigne le nombre de places de $K$ au-dessus de $\ell$.\medskip

Le but de cette note est de préciser la valeur de $\wi e_K$ indépendamment de toute conjecture ainsi que les cas d'égalité entre le groupe $\;\wE_K$ et son sous-module {\em naïf} \;$\Zl\otimes_{\ZZ}\wi E_K$.

\newpage
\noindent{\large \bf 1.Bref rappel  sur les unités logarithmiques}
\medskip

Soient $\ell$ un nombre premier donné et $K$ un corps de nombres.  À chaque place finie  $\p$ de $K$, il est attaché dans \cite{J28} une application à valeurs dans $\Zl$, définie sur le groupe multiplicatif $K_\p^\times$ du complété de $K$ en $\p$ par la formule:

\centerline{$\tilde\nu_\p (x_\p)\, =\,\nu_\p (x_\p)$, pour $\p\nmid\ell$;\quad et \quad $\tilde\nu_\p (x_\p)\, =  -\frac{\Log_\ell\, N_{K_\p/\QQ_p}(x_\p)}{\deg\, \p}$, pour $\p|\ell$;}\smallskip

\noindent où $\Log_\ell$ désigne le logarithme d'Iwasawa et $\deg\p$ est un facteur de normalisation, dont l'expression exacte est sans importance ici, destiné à assurer que l'image de $K_\p^\times$ soit dense dans $\Zl$. Cette application induit un morphisme surjectif du compactifié $\ell$-adique $\R_{K_\p}=\varprojlim  K_\p^\times/K_\p^{\times\ell^n}$ de  $K_\p^\times$ sur $\Zl$ dont le noyau, dit {\em sous-groupe des unités logarithmiques} de $\R_{K_\p}$,\smallskip

\centerline{$\wU_{K_\p}\,=\,\{u_\p\in\R_{K_\p}\,|\, \tilde\nu_\p(u_\p)=0\}$}\smallskip

\noindent  s'identifie par la Théorie $\ell$-adique locale du corps de classes (cf. \cite{J31}) au sous groupe normique de $\R_{K_\p}$ associé à la $\Zl$-extension cyclotomique $K_\p^c$ de $K_\p$.\smallskip

Soit maintenant $\J _K$ le {\em $\ell$-adifié  du groupe  des idèles de $K$}, i.e. le produit $\J_K=\prod_\p ^{\rm res}\R_{K_\p}$ des compactifés  $\R _{K_\p}$ des groupes multiplicatifs des complétés $K_\p$, restreint aux familles $(x_\p)_\p$ dont presque tous les éléments tombent dans le sous-groupe unité $\,\wU_K=\prod_\p \, \wU _{K_\p}$.
La Théorie $\ell$-adique globale du corps de classes établit un isomorphisme de groupes topologiques compacts entre le $\ell$-groupe des classes d'idèles $\,\C_K$ défini comme quotient\smallskip

\centerline{$\C_K=\J_K/\R_K$}\smallskip

\noindent de $\J_K$ par son sous-groupe principal $\R_K=\Zl\otimes_\ZZ K^\times$ et le groupe de Galois $G_K=\Gal(K^{ab}/K)$ de la pro-$\ell$-extension abélienne maximale de $K$. Dans la correspondance ainsi établie (cf. \cite{J28,J31}):\smallskip

\begin{itemize}
\item[(i)] Le groupe de normes associé à la $\Zl$-extension cyclotomique $K^c$ du corps $K$ est le sous-groupe des idèles de degré nul: $\wJ_K=\{\x=(x_\p)_\p\in\J_K\,|\, \deg(\x)=\sum_\p\wi\nu_\p(x_\p)\deg\,\p=0\}$.
\item[(ii)] Le groupe de normes associé à la plus grande sous-extension $K^{lc}$ de $K^{ab}$ qui est localement cyclotomique (i.e. complètement décomposée sur $K^c$ en chacune de ses places) est le produit $\,\wU_K\R_K$ du sous-groupe $\,\wU_K=\prod_\p \, \wU _{K_\p}$ des unités logarithmiques locales et de $\R_K$.
\item[(iii)] En particulier, le groupe de Galois $\Gal(K^{lc}/K^c)$ s'identifie au quotient $\,\wCl_K= \wJ_K/\wU_K\R_K$, lequel peut être regardé comme quotient du groupe $\wDl _K = \wJ _K /\wU_K$ des diviseurs logarithmiques de degré nul par son sous-groupe principal $\Pl_K=\R_K\wU_K/\wU_K$, le numérateur $\,\wDl _K$  s'identifiant au sous-groupe $\wi\oplus_\p \,\Zl\,\p$ des diviseurs de degré nul de la somme formelle $\oplus_\p \,\Zl\,\p$.
\item[(iv)] Et le noyau $\wE _K=\R_K \cap\, \wU_K$ du morphisme $\wi\div:\;x\mapsto \sum_\p\wi\nu_\p(x)\,\p$ de $\R_K$ dans $\wDl_F$ est le sous-groupe des normes cyclotomiques globales de $\R_K$.
\end{itemize}
\smallskip

\noindent Nous disons que $\,\wCl_K$ est le {\em $\ell$-groupe des classes logarithmiques} du corps $K$ et que $\,\wE_K$ est le {\em $\ell$-groupe des unités logarithmiques globales}.\smallskip

Dans ce contexte, la conjecture de Gross-Kuz'min (cf. \cite{J55} \S2) se présente comme suit:

\begin{ThC*}\label{Ulog}
Le groupe $\,\wE_K$ des unités logarithmiques de $K$ est le produit:

\centerline{ $\wE_K= \mu^{_{(\ell)}}_K  \Zl^{r_K+c_K+\delta_K^{\mathscr G}}$}\smallskip

\noindent du $\ell$-sous-groupe $ \mu^{_{(\ell)}}_K$ des racines globales de l'unité et d'un $\Zl$-module libre de rang $r_K+c_K+\delta_K^{\mathscr G}$, où  $r_K$ et $c_K$ sont respectivement les nombres de places réelles et complexes de $K$; et $\delta_K^{\mathscr G}=\dim^{\phantom{c}}_{\Zl} \wCl_K$  mesure le défaut dans $K$ de la conjecture de Gross-Kuz'min  pour  le premier $\ell$, laquelle postule précisément la finitude du groupe $\wCl_K$ ou, de façon équivalente, l'égalité: $\dim\wE_K=r_K+c_K$.
\end{ThC*}

Il est bien connu (cf. e.g. \cite{FG,J10,J17,Kuz}) que la conjecture de Gross-Kuz'min est satisfaite pour tous les premiers $\ell$ dès lors que $K$ est abélien, en vertu du théorème d'indépendance de Baker-Brumer. Mais c'est aussi le cas pour des raisons triviales (du fait de la formule du produit), lorsque le corps $K$ admet une seule place au-dessus du premier $\ell$:  la valuation logarithmique $\tilde\nu_\p$ définie plus haut coïncide, en effet, avec la valuation $\nu_\p$ au sens ordinaire dès que $\p$ ne divise pas $\ell$. Il suit de là que le groupe $\,\wE_K$ est contenu dans le tensorisé $\,\E'_K=\Zl\otimes_{\ZZ}E'_K$ du groupe des $\ell$-unités de $K$. D'où leur égalité pour $\dim_\ZZ E'_K=r_K+c_K$, i.e. lorsque $K$ admet $l_K=1$ places au-dessus de $\ell$.\par
Dans la pratique, le calcul du groupe $\,\wCl_K$ est parfaitement effectif (cf. \cite{BJ,DJ+}).

\newpage
\noindent{\large \bf2.  Normes cyclotomiques naïves}
\medskip

Intéressons-nous maintenant au groupe $\tilde E_K$ des $\ell$-normes cyclotomiques au sens naïf, c'est-à-dire à l'intersection

\centerline{$\wi E_K =\bigcap_{n\in\NN} N_{K_n/K}(K_n^\times)$}\smallskip

\noindent des sous-goupes normiques de $K^\times$ attachés aux étages finis de la $\Zl$-extension cyclotomique $K^c/K$.\smallskip

Par le principe de Hasse, la condition normique se lit localement, de sorte que $\wi E_K$ est encore le noyau dans $K^\times$ des valuations logarithmiques $\tilde\nu_\p$ présentées plus haut:\smallskip

\centerline{$\wi E_K\,=\,\{x\in K^\times \,|\, \tilde\nu_\p(x)=0 \quad\forall \p\nmid\infty\}\,=\,\{\varepsilon\in E'_K \,|\, \tilde\nu_\p(\varepsilon)=0 \quad\forall \p\mid\ell\}$.}\smallskip

En particulier $\wi E_K$ est un sous-groupe du groupe $E'_K$ des $\ell$-unités de $K$ qui contient les racines de l'unité ainsi que les puissances de $\ell$, ce qui se traduit par les deux inclusions:\smallskip

\centerline{$\mu_K^{\phantom{lc}}\ell^\ZZ \, \subset \, \wi E_K \,\subset\, E'_K \simeq \mu_K^{\phantom{lc}}\ZZ^{r_K+c_K+l_K-1}$,}\smallskip

\noindent où $r_K$, $c_K$ et $l_K$ désignent respectivement les nombres de places réelles, complexes et $\ell$-adiques de $K$, de sorte que l'on a:

\centerline{$\wi E_K \simeq  \mu_K^{\phantom{lc)}}\ZZ^{\tilde e_K}$}\smallskip

\noindent et la dimension $\wi e_K$ de $\wi E_K$ est comprise entre 1 et $r_K+c_K+l_K-1$.\smallskip

Observons tout de suite que les deux bornes sont atteintes sous des hypothèses convenables:

\begin{Prop} Soit $\wi e_K=\dim_\ZZ\wi E_K$ la dimension du groupe des $\ell$-normes cyclotomiques naïves.
\begin{itemize}
\item[(i)] Si le corps $K$ possède une unique place au-dessus de $\ell$, $\wi E_K$ coïncide avec le groupe des $\ell$-unités $E'_K$ et il vient: $\wi e_K=r_K+c_K+l_K-1=r_K+c_K$.
\item[(ii)] Si, tout au contraire, la place $\ell$ est complètement décomposée dans $K/\QQ$, on a: $\wi E_K=\mu_K^{\phantom{lc}}\ell^\ZZ$ et $\wi e_K=1$. Nous disons alors que $K$ est totalement $\ell$-adique.
\end{itemize}
\end{Prop}

\Preuve Le cas (i) ayant été traité à la fin de la section précédente, concentrons nous sur le cas (ii). Par hypothèse, la norme locale $N_{K_\l/\Ql}$ aux places $\l$ au-dessus de $\ell$ étant triviale, les valuations logarithmiques correspondantes $\wi\nu_\l$ sont proportionnelles au logarithme d'Iwasawa $\Log_\ell$ et il suit:\smallskip

\centerline{$\wi E_K= \{\varepsilon\in E'_K \,|\, Log_\ell(\varepsilon)=0\} = \mu_K^{\phantom{lc}}\ell^\ZZ$,}\smallskip

\noindent puisque, du fait de l'hypothèse de complète décomposition, les seules puissances fractionnaires de $\ell$ contenues dans $K$ sont les puissances entières.\medskip

Plus généralement, en présence de conjugaison $\ell$-adique (cf. \cite{J47} pour cette notion), il vient:

\begin{DTh}\label{CLA}
Convenons de dire qu'un corps de nombres $K$ est à conjugaison $\ell$-adique lorsque c'est une extension non-décomposée aux places au-dessus de $\ell$ d'un sous-corps $k$ complètement décomposé au-dessus de $\ell$. Pour un tel corps $K$ la dimension de $\wi E_K$ est donnée par:\smallskip

\centerline{$\wi e_K= (r_K+c_K)-(r_k+c_k-1)$,}\smallskip

\noindent où $r_K$, $r_k$, $c_K$ et $c_k$ sont les nombres respectifs de plongements réels et complexes de $K$ et $k$.
\end{DTh}

\Preuve Notons $V'_K=\QQ\otimes_\ZZ E'_K$ et $V'_k=\QQ\otimes_\ZZ E'_k$ les $\QQ$-espaces vectoriels construits sur les groupes de $\ell$-unités respectifs de $K$ et $k$, puis $V^*_{K/k}$ le noyau de l'opérateur norme $\nu=N_{K/k}$ de $V'_K$ sur $V'_k$.

Par construction, le sous-espace $\wi V_K=\QQ\otimes_\ZZ \wi E_K$ de $V'_K$ construit sur les unités logarithmiques est donné par:

\centerline{$\wi V_K = \{\varepsilon \in V'_K \, |\, N_{K/k}(\varepsilon) \in \ell^\QQ\} $,}\smallskip

\noindent comme préimage par la norme $N_{K/k}$ du noyau du logarithme. Il vient donc:\smallskip 

\centerline{$\dim_\ZZ\wi E_K = \dim_\QQ\wi V_K=1+\dim_\QQ V^*_{K/k}= 1+(r_K+c_K+l_K-1)-(r_k+c_k+l_k-1)$;}\smallskip

\noindent d'où le résultat annoncé, puisque $K$ et $k$ ont même nombre $l_K=l_k$ de places au-dessus de $\ell$.

\begin{Sco}
Sous les hypothèses du Théorème \ref{CLA}, on a l'égalité $\,\wi e_K=r_K+c_K$ (et donc, sous la conjecture de Gross-Kuz'min, l'identité $\,\wE_K=\Zl\otimes_\ZZ\wi E_K$) si et seulement si le sous-corps totalement $\ell$-adique $k$ est ou bien le corps des rationnels $\QQ$ ou bien un corps quadratique imaginaire $\QQ[\sqrt{-d}]$.
\end{Sco}

\newpage
\noindent{\large \bf3. Étude du cas galoisien}
\medskip

Examinons maintenant le cas où le corps de nombres considéré est galoisien de degré $n$ sur $\QQ$ et notons $G=\Gal(K/\QQ)$ son groupe de Galois.

Faisons choix d'une des places à l'infini, disons $\p_{\!_\infty\!}$ de $K$ et notons $D_{\p_{\!_\infty\!}}$ son groupe de décomposition; notons de même $\p_{\!_\ell}$ l'une des places de K au-dessus de $\ell$ et $D_{\p_{\!_\ell}}$ son groupe de décomposition.
Notons enfin $d_{\infty}$ l'ordre de  $D_{\p_{\!_\infty\!}}$ (qui vaut 1 ou 2) et $d_\ell$ celui de $D_{\p_{\!_\ell}}$.\smallskip

D'après le Théorème de Herbrand, le caractère du $\QQ[G]$-module $V'_K=\QQ\otimes E'_K$ construit sur le groupe des $\ell$-unités de $K$ est donné par:

\centerline{$\chi_{E'_K}^{\phantom{lc}}=\chi_\infty^{\phantom{lc}}+\chi_\ell^{\phantom{lc}}-1$,}\smallskip

\noindent où $\chi_\infty^{\phantom{lc}}=\Ind_{D_{\p_{\!_\infty\!}}}^G 1_{D_{\p_{\!_\infty\!}}}$ est l'induit à $G$ du caractère de la représentation unité du sous-groupe de décomposition de $\p_{\!_\infty\!}$; et $\chi_\ell^{\phantom{lc}}=\Ind_{D_{\p_{\!_\ell}}}^G 1_{D_{\p_{\!_\ell}}}$ son analogue pour $\p_{\!_\ell}$; enfin 1 désigne le caractère unité.

L'objet de cette section est de calculer le caractère $\chi_{\wi E_K}$ du $\QQ[G]$-module $\wi V_K=\QQ\otimes_\ZZ \wi E_K$.\medskip

Or, d'un côté les 1-composantes isotypiques respectives de $V'_K$ et de $\wi V_K$ s'identifient au sous-espace $V'_\QQ = \wi V_\QQ = \ell^{\,\QQ}$; de sorte que $\chi_{E'_K}^{\phantom{lc}}$ et $\chi_{\wi E_K}^{\phantom{lc}}$ contiennent une seule fois le caractère unité.\smallskip

Et d'un autre côté, si $E_K^*$ désigne le noyau dans $E'_K$ de la norme $N_{K/\QQ}$, le $\QQ$-espace associé $V^*_K=\QQ\otimes_\ZZ E^*_K$ n'est autre que le sous-$\QQ[G]$-module de $V'_K$ annulé par l'idempotent $e_1^{\phantom{lc}}=\frac{1}{n}\sum_{\gamma\in G}\gamma$; et le caractère de $V^*_K$ est donné par: 

\centerline{$\chi_{E^*_K}^{\phantom{lc}}=(\chi_\infty^{\phantom{lc}}-1)+(\chi_\ell^{\phantom{lc}}-1)$.}\smallskip

\noindent L'expression des valeurs absolues logarithmiques donnée au début de la section 1 montre alors que le sous-espace $\wi V^*_K=\QQ\otimes_\ZZ \wi E^*_K$ de $V^*_K$ construit sur les unités logarithmiques de norme 1 est l'intersection des noyaux dans $V^*_K$ de l'idempotent associé au sous-groupe $D_{\p_{\!_\ell}}$ \smallskip

\centerline{$e_{\p_{\!_\ell}}\;=\;\frac{1}{n} \underset{\sigma\in G}{\sum}\,\chi_\ell^{\phantom{c}}(\sigma^{-1})\sigma\;=\;\frac{1}{d_\ell} \underset{\tau\in D_{\p_{_\ell}}}{\sum}\,\tau$}

\noindent et de ses conjugués $e_{\p^\gamma_{\!_\ell}}=\gamma e_{\p_{\!_\ell}}\gamma^{-1}$, lorsque $\gamma$ décrit $G$. En résumé:

\begin{Prop}
Le $\QQ[G]$-module $\wi V_K=\QQ\otimes_\ZZ \wi E_K$ construit sur les unités logarithmiques naïves est la somme directe  $\wi V_K =\wi V_\QQ \oplus \wi V^*_K$, où
\begin{itemize}
\item $\wi V_\QQ = \ell^\QQ$ est l'image de $\wi V_K$ par l'idempotent central $e_1^{\phantom{lc}}=\frac{1}{n}\sum_{\gamma\in G}\gamma$ attaché à la norme;
\item $\wi V_K^*$ est le plus grand sous-module de $V^*_K$ qui est tué par l'idempotent $e_{\p_{\!_\ell}}\;=\;\frac{1}{d_\ell} \underset{\tau\in D_{\p_{_\ell}}}{\sum}\,\tau$.
\end{itemize}
\end{Prop}

À partir de là, la détermination du caractère $\chi_{\wi E_K}$ est un pur problème de théorie des représentations où l'arithmétique n'a plus aucune part: on dispose d'un groupe fini $G$, d'un $\QQ[G]$ module projectif  $X=V^*_K$ de caractère $\chi=\chi_{E^*_K}^{\phantom{lc}}=(\chi_\infty^{\phantom{lc}}-1)+(\chi_\ell^{\phantom{lc}}-1)$, d'un sous-groupe $H=D_{\p_{\!_\ell}}$ de $G$; et on cherche le caractère du plus grand sous-module de $X$ qui est annulé par l'idempotent $e_H=\frac{1}{|H|}\sum_{\tau\in H}\,\tau$ de l'agèbre $\QQ[G]$ construit sur les éléments de $H$.\smallskip

Écrivons donc $\chi=\sum_{i\in I} n_i\chi_i$ la factorisation irréductible de $\chi$ et $X\simeq\bigoplus_{i\in I}\bigoplus_{j=1}^{n_i} X_i$ comme somme directe d'idéaux à gauche minimaux deux à deux non isomorphes. Pour un $i$ donné, on a:
\begin{itemize}
\item ou bien $e_H X_i=0$; et $X_i$ est annulé par $e_H$ (et tous ses conjugués);
\item ou bien $e_H X_i \ne 0$; et le seul sous-module de $X_i$ qui est annulé par $e_H$ est le module nul.
\end{itemize}\smallskip

En fin de compte, le sous-module cherché $\wi X$ est donc donné par: $\wi X \simeq \bigoplus_{e_HX_i=0}\bigoplus_{j=1}^{n_i} X_i$.
Notons que la condition $e_HX_i=0$ se lit encore $e_H A_i=0$, si $A_i=\QQ[G]e_{\chi_i}$ est le facteur simple correspondant au caractère $\chi_i$, autrement dit: $e_He_{\chi_i}=0$. En pratique, on obtient donc $\wi V^*_K$ à partir de $V^*_K$ en éliminant les composantes irréductibles représentés dans le $\QQ[G]$-module $\QQ[G/H]\simeq\QQ[G]e_H$ construit sur les classes à gauche modulo $H$. Ainsi:

\begin{Th}
Pour $K/\QQ$ galoisienne le caractère du groupe des normes cyclotomiques naïves est:\smallskip

\centerline{$\chi_{\wi E_K}^{\phantom{lc}}= 1+(\chi_{\infty\!}^{\phantom{lc}}\wedge\bar\chi_\ell^{\phantom{lc}})$,}\smallskip

\noindent où $\bar \chi_\ell^{\phantom{c}}$ désigne la partie du caractère régulier $\chi_{\rm r\acute{e}g}$ qui est étrangère à $\chi_\ell^{\phantom{c}}$.
\end{Th}

\newpage
\noindent{\large \bf4. Discussion des cas d'égalité: $\,\wE_K=\Zl\otimes_\ZZ \wi E_K$}
\medskip

Le groupe des normes cyclotomique naïves $\wi E_K$ est un sous-module pur du groupe des $\ell$-unités $E'_K$, comme noyau des valuations logarithmiques $\wi\nu_\p$ attchées aux places $\p|\ell$; son $\ell$-adifié $\Zl\otimes_\ZZ \wi E_K$ est donc un sous-module pur du $\ell$-adifié $\E'_K=\Zl\otimes_\ZZ E'_K$ de $E'_K$. Comme il est contenu dans $\,\wE_K$, c'est en particulier un sous-module pur de $\,\wE_K$ et l'égalité $\Zl\otimes_\ZZ \wi E_K=\wE_K$ a lieu si et seulement si les deux modules ont le même rang, i.e. lorsque le $\ZZ$-rang de $\wi E_K$ coïncide avec le $\Zl$-rang de $\,\wE_K$.

Dans le cas galoisien, cela revient à dire que les deux groupes définissent le même caractère. Or, celui du groupe des unités logarithmiques est calculé dans \cite{J28}: il est égal à $\chi_\infty^{\phantom{lc}}$ augmenté du caractère de défaut de la conjecture de Gross-Kuz'min. Il vient donc directement:

\begin{Th}
Dans une extension galoisienne $K$ de $\QQ$, l'égalité  $\,\wE_K=\Zl\otimes_\ZZ \wi E_K$ a lieu si et seulement si les deux conditions suivantes sont réalisées:
\begin{itemize}
\item[(i)] le corps $K$ vérifie la conjecture de Gross-Kuz'min (pour le premier $\ell$);
\item[(ii)] on a: $\chi_\infty^{\phantom{c}}\wedge\chi_\ell^{\phantom{c}}=1$.
\end{itemize}
\end{Th}

\begin{Cor}
Pour $K$ galoisien réel, l'égalité  $\,\wE_K=\Zl\otimes_\ZZ \wi E_K$ a lieu si et seulement si le corps $K$ possède une seule place au-dessus de $\ell$; auquel cas, il vient: $\wi E_K=E'_K$ et $\,\wE_K=\Zl\otimes_\ZZ \wi E_K=\E'_K$.
\end{Cor}

\Preuve Si $K$ est réel, on a $\chi_\infty^{\phantom{c}}=\chi_{\rm r\acute{e}g}$ et l'assertion (ii) s'écrit: $\chi_\ell^{\phantom{c}}=1$; ce qui se traduit par le fait que $K$ admet une unique place au-dessus de $\ell$. D'où le résultat en vertu de la Proposition 1. Notons que, dans ce contexte, la condition (i) est alors automatiquement satisfaite.

\begin{Cor}
Si $K$ est à conjugaison complexe, i.e. si c'est une extension quadratique totalement imaginaire d'un sous-corps $K_{\infty\!}$ totalement réel, 
l'égalité  $\,\wE_K=\Zl\otimes_\ZZ \wi E_K$ a lieu si et seulement si
\begin{itemize}
\item[(i)] ou bien $K$ ne possède qu'une seule place au-dessus de $\ell$;
\item[(ii)] ou bien $K$ est composé direct d'un sous-corps galoisien réel $K_{\infty\!}$ qui ne possède qu'une place au-dessus de $\ell$ et d'un sous-corps quadratique imaginaire $k$ qui en possède exactement 2.
\end{itemize}
\end{Cor}

\noindent Avant d'établir ce résultat, observons que nous avons, par un argument de points fixes immédiat:

\begin{Lem}
Si l'égalité $\,\wE_K=\Zl\otimes_\ZZ \wi E_K$ a lieu pour un corps de nombres quelconque donné $K$, elle vaut aussi pour chacun de ses sous-corps.
\end{Lem}

\PreuveCor Si $K$ est à conjugaison complexe, le sous-groupe de décomposition $D_{\infty\!}$ des places à l'infini est normal et d'ordre 2 et son sous-corps des points fixes $K_{\infty\!}$ est alors un corps galoisien réel auquel nous pouvons appliquer le corollaire précédent. Si donc $K$ satisfait l'égalité, $K_{\infty\!}$ la satisfait aussi et ne possède qu'une seule place au-dessus de $\ell$. Ainsi pour chaque place $\ell$-adique de $K$, le sous-groupe de décomposition s'envoie surjectivement sur le quotient $G/D_\infty$. Il est donc d'indice 1 ou 2 dans $G$ et, par conséquent, normal. Notons le $D_\ell$.
\begin{itemize}
\item Dans le premier cas ($D_\ell=G$), le corps  $K$ ne possède qu'une seule place au-dessus de $\ell$; et il vient, comme plus haut:  $\wi E_K=E'_K$ et $\,\wE_K=\Zl\otimes_\ZZ \wi E_K=\E'_K$.
\item Dans le second cas  ($D_\ell\ne G$), le groupe $G$ est alors le produit direct des sous-groupes  $D_{\infty\!}$ et $D_\ell$, et $K$ le composé direct de leurs sous-corps invariants respectifs $K_\infty$ et $k$. Ici encore le corps réel $K_\infty$ vérifie les hypothèses du corollaire précédent, de sorte qu'on a:  $\wi E_{K_\infty\!}=E'_{K_\infty}$ et $\,\wE_{K_\infty\!}=\Zl\otimes_\ZZ \wi E_{K_\infty\!}=\E'_{K_\infty}$. Comme, sous la conjecture de Gross-Kuz'min, les groupes $\,\wE_{K_\infty}$ et $\,\wE_K$ ont même $\Zl$-rang $r_{K_\infty\!}=c_K$, il suit bien: $\,\wE_K=\mu_K^{_{(\ell)}}\wE_{K_\infty}=\Zl\otimes_\ZZ \wi E_K$.
\end{itemize}
En fin de compte, il reste seulement à vérifier dans ce dernier cas que le corps $K=k K_\infty$ vérifie la conjecture de Gross-Kuz'min. Or, cela résulte de la proriété générale suivante:

\begin{Sco}
Soit $K/k$ une extension (non nécessairement galoisienne) de corps de nombres, non décomposée aux places au-dessus de $\ell$. Alors la conjecture de Gross-Kuz'min pour le premier $\ell$ est vérifiée dans $K$ dès qu'elle l'est dans le sous-corps $k$ (par exemple dès que $k$ est abélien sur $\QQ$).
\end{Sco}

\Preuve Elle est très simple: le corps $K$ et son sous-corps $k$ ayant même nombre de places au-dessus de $\ell$, le groupe des diviseurs logarithmiques de $k$ construits sur les places au-dessus de $\ell$ est d'indice fini dans son homologue de $K$. Ainsi $\wCl_{k}$ et $\wCl_K$ sont-ils simultanément finis ou pas.

\newpage
\noindent{\sc Index des principales notations}\medskip
\medskip

$\ell$: un nombre premier quelconque;\

$K$  : un corps de nombres arbitraire; $K_\p$: le complété de $K$ en la place $\p$;\
	
$K^c = \cup_{n \in \mathbb N} K_n$ avec $[K_n :K]= \ell^n$: la $\Zl$-extension cyclotomique de $K$;\

$E_K$ : le groupe des unités de $K$ et $E'_K$ le groupe des $\ell$-unités;\

$\wi E_K =\bigcap_{n\in\NN} N_{K_n/K}(K_n^\times)$: le groupe des normes cyclotomiques naïves de $K$;\

$\mu_K^{\phantom{c}}$: le groupe des racines de l'unité de $K$ et $\mu_K^{_{(\ell)}}$ son $\ell$-sous-groupe de Sylow;\

$r_K, c_K, l_K$: les nombres respectifs de places réelles, complexes et $\ell$-adiques du corps $K$;\

$\R_K=\Zl\otimes_\ZZ K^\times$: le $\ell$-adifié du groupe multiplicatif du corps $K$;\

$\R_{K_\p} = \varprojlim K^\times_\p/K^{\times \ell^m}_\p\!\!$: le compactifié $\ell$-adique du groupe multiplicatif $K^\times_\p$;\

$\U_{K_\p}$: le sous-groupe unité et $\wU_{K_\p}$ le groupe des normes cyclotomiques dans $\R_{K_\p}$;\

$\J_K = \prod^\mathrm{res}_\p \R_{K_\p}$ : le $\ell$-adifié du groupe des idèles de $K$;\

$\U_K=\prod_\p \U_{K_\p}$: le sous-groupe unité de $\J_K$;\

$\wU_K=\prod_\p \wU_{K_\p}$: le sous-groupe des normes cyclotomiques dans $\J_K$;\

$\wCl_K=\J_K/\wU_K\R_K$: le $\ell$-groupe des classes logarithmiques du corps $K$;

$\E'_K=\Zl\otimes_\ZZ E'_K$: le $\ell$-adifié du groupe des $\ell$-unités de $K$;\

$\E_K=\Zl\otimes_\ZZ E_K=\R_K\cap\,\U_K$: le $\ell$-adifié du groupe des unités de $K$;\

$\wE_K=\R_K\cap\,\wU_K$: le $\ell$-groupe des unités logarithmiques de $K$;\

$\delta^{\mathscr G}_K=\dim_{\Zl}\wCl_K$: le défaut de la conjecture de Gross-Kuz'min dans $K$.\



\def\refname{\normalsize{\sc  Références}}

\bigskip\noindent
{\small
\begin{tabular}{l}
{Jean-François {\sc Jaulent}}\\
Institut de Mathématiques de Bordeaux \\
Université de {\sc Bordeaux} \& CNRS \\
351, cours de la libération\\
F-33405 {\sc Talence} Cedex\\
courriel : Jean-Francois.Jaulent@math.u-bordeaux1.fr 
\end{tabular}
}

 \end{document}